\documentclass[12pt]{amsart}
\newcommand{\diag}{\mathop{\rm diag\,}\nolimits}
\newcommand{\f}{\mathbf f}
\newcommand{\g}{\mathbf g}

\begin{document}
\title{A formula for inverse matrix}
\author{Alexander Roi Stoyanovsky}
\begin{abstract}
We give a formula for the inverse matrix to an infinite matrix with possibly noncommutative entries, generalizing the  Newton 
interpolation formula and the Taylor formula. 
\end{abstract}
\maketitle

\section*{Introduction}

Many important problems of classical analysis can be interpreted as the problems of inverting a linear operator between infinite dimensional vector spaces.
This viewpoint has been emphasized already in the classical book [1]. If we choose infinite bases of the vector spaces in question, the problem is to solve
an infinite system of linear equations with infinitely many variables, or to invert a matrix of infinite order.

For example, let $x_0, x_1, \ldots$ be pairwise distinct numbers, and let $f(x)=\sum\limits_{k=0}^\infty a_kx^k$ be an analytical function. The
interpolation problem, i.~e. the problem of finding the function $f(x)$ from the equalities $f(x_i)=f_i$, $i=0,1,\ldots$, is equivalent to the problem of inverting the 
infinite matrix $(x_i^k)_{i,k=0}^\infty$.

The answer to the interpolation problem is given by the classical Newton interpolation formula. Define the {\it difference derivatives} (or simply {\it differences})
$\Delta^m f(x_0,\ldots,x_m)$ by the formulas
\begin{equation}
\begin{aligned}{}
&\Delta f(x_0,x_1)=\frac{f(x_0)-f(x_1)}{x_0-x_1},\\
&\Delta^m f(x_0,\ldots,x_m)=\frac{\Delta^{m-1}f(x_0,\ldots,x_{m-1})-\Delta^{m-1}f(x_0,\ldots,x_{m-2},x_m)}{x_{m-1}-x_m}.
\end{aligned}
\end{equation}
Then the Newton formula (for a modern account, see e.~g. [2]) is
\begin{equation}
\begin{aligned}{}
&f(x)=\sum\limits_{m=0}^n (x-x_0)\ldots(x-x_{m-1})\Delta^m f(x_0,\ldots,x_m)+R_n(x),\\
&R_n(x_0)=\ldots=R_n(x_n)=0.
\end{aligned}
\end{equation}

A remarkable and somewhat unexpected property of the differences $\Delta^m f(x_0,\ldots,x_m)$ is that they are symmetric functions of $x_0,\ldots,x_m$. The present work arose 
from an attempt to understand and find a natural context for this property. The result is a formula for the inverse matrix generalizing the Newton formula.

Another remarkable property of the difference derivatives 
$$\Delta^m f(x_0,\ldots,x_m)$$ 
is that they have finite limits as the points $x_0,\ldots,x_m$ merge.
For example, when $x_1,\ldots,x_m\to x_0$, one has
\begin{equation}
\Delta^m f(x_0,\ldots,x_m)\to\frac{f^{(m)}(x_0)}{m!},
\end{equation}
where $f^{(m)}(x)$ is the $m$-th derivative of $f$, and the Newton formula tends to the Taylor formula
\begin{equation}
\begin{aligned}{}
&f(x)=\sum\limits_{m=0}^n (x-x_0)^m\frac{f^{(m)}(x_0)}{m!}+R_n(x),\\
&R_n(x_0)=R_n'(x_0)=\ldots=R_n^{(n)}(x_0)=0.
\end{aligned}
\end{equation}
The Taylor formula is also a particular case of our formula. 

It turns out that the most simple and natural formulas arise if we allow the product operation to be noncommutative, and use the theory of quasideterminants due to I.~M.~Gelfand, V.~S.~Retakh et al. [3]. 

The paper is organized as follows. In \S1 we derive the generalized Newton formula. In \S2 we give examples and applications of this formula, including the Newton formula and the Taylor formula.

\section{Generalization of the Newton formula}

\subsection{Notations} Let $V$ be a left (topological) vector space over a skew-field with a basis $v^k, k=0,1,\ldots$, let $V'$ be the dual right (topological) 
vector space with a basis $p_i, i=0,1,\ldots$, let $\f\in V$ be a vector, and let $\g\in V'$ be a covector. 
Denote $y_i^k=(v^k,p_i)$, $f_i=(\f,p_i)$, $g^k=(v^k,\g)$.

Denote by $D$ the infinite matrix $(y_i^k)$,  $i, k=0,1,2,\ldots$, by $f$ the row vector $(f_i)$, $i=0,1,\ldots$, and by $g$ the column vector $(g^k)$, $k=0,1,\ldots$.
Denote by $y^k=(y_0^k,y_1^k,\ldots)$ the $k$-th row of $D$, and by $y_i=(y_i^0,y_i^1,\ldots)$ the $i$-th column of $D$. 

Denote by $D_{i_0\ldots i_m}^{k_0\ldots k_m}$ the $(m+1)\times(m+1)$-submatrix of $D$ with the columns $i_0,\ldots,i_m$ and the rows $k_0,\ldots,k_m$. 

Denote by $|M|_i^j$ the $(i,j)$-quasideterminant of an $n\times n$-matrix $M$ [3]. In the commutative case, one has $|M|_i^j=(-1)^{i+j}\det M/\det M_i^j$, where $M_i^j$ is the $(n-1)\times(n-1)$-submatrix
obtained from $M$ by removing the $i$-th column and the $j$-th row.

\subsection{The biorthogonalization process}
{\bf Theorem 1.} {\it For a generic matrix $D$, there exist the unique infinite upper triangular matrix $A=(a_m^i)$ and the unique infinite lower triangular matrix $C=(c_k^m)$ such that the vectors 
\begin{equation}
w^m=\sum\limits_{k=0}^m c_k^m v^k\in V,\ \  q_m=\sum\limits_{i=0}^m p_i a_m^i\in V'
\end{equation}
form biorthogonal bases in $V$ and $V'$, i.~e. satisfy the properties  }
\begin{equation}
\begin{aligned}{}
(w^m,q_{m'})&=0\text{ for }m\ne m',\\
(w^m,p_m)&=(v^m,q_m)=1.
\end{aligned}
\end{equation}

{\it Proof.} The vectors $w^m$, $q_m$ are the unique linear combinations of respectively $v^0,\ldots,v^m$ and $p_0,\ldots,p_m$ satisfying the properties
\begin{equation}
\begin{aligned}{}
&(w^m,p_i)=0,\ \ 0\le i\le m-1,\ \ (w^m,p_m)=1;\\
&(v^k,q_m)=0,\ \ 0\le k\le m-1, \ \ (v^m,q_m)=1.\ \ \qed
\end{aligned}
\end{equation}

\subsection{Difference derivatives}
Let $i_0, i_1, \ldots$ and $k_0, k_1, \ldots$ be two permutations of the sequence $0,1,\ldots$. 
Consider the basis $\widetilde p_0=p_{i_0},\widetilde p_1=p_{i_1},\ldots$ in $V'$ and the basis $\widetilde v^0=v^{k_0},\widetilde v^1=v^{k_1},\ldots$ in $V$.
Let us apply the biorthogonalization process, i.~e. Theorem 1, to these bases. We obtain the infinite matrices $A=(a_m^i)$, $C=(c_k^m)$, and the biorthogonal bases $w^m$, $q_m$. 

{\bf Definition.} Let us call the numbers 
\begin{equation}
\begin{aligned}{}
&(\Delta_R^m f)_{i_0\ldots i_m}^{k_0\ldots k_m}=\sum\limits_{j=0}^m f_{i_j}a_m^j=(\f,q_m),\\ 
&(\Delta_L^m g)_{i_0\ldots i_m}^{k_0\ldots k_m}=\sum\limits_{l=0}^m c_l^m g^{k_l}=(w^m,\g) 
\end{aligned}
\end{equation}
by respectively the {\it $m$-th right difference derivative} of the vector $f$ and the {\it $m$-th left difference derivative} of the vector $g$ (corresponding to the sequences $i_0,i_1,\ldots$ and 
$k_0,k_1,\ldots$).
\medskip

{\it Recurrence formulas for difference derivatives.}
The following Theorem justifies the term ``difference derivatives''.

{\bf Theorem 2.} a) {\it We have}
\begin{equation}
\begin{aligned}{}
&(\Delta_R^0 f)_i^k=f_i(y_i^k)^{-1},\\
&(\Delta_R^1 f)_{ii'}^{kk'}=((\Delta_R^0 f)_i^k-(\Delta_R^0 f)_{i'}^k)((\Delta_R^0 y^{k'})_i^k-(\Delta_R^0 y^{k'})_{i'}^k)^{-1},\\
&(\Delta_R^m f)_{i_0\ldots i_m}^{k_0\ldots k_m}=((\Delta_R^{m-1}f)_{i_0\ldots i_{m-1}}^{k_0\ldots k_{m-1}}-(\Delta_R^{m-1}f)_{i_0\ldots i_{m-2}i_m}^{k_0\ldots k_{m-1}})\\
&\times((\Delta_R^{m-1}y^{k_m})_{i_0\ldots i_{m-1}}^{k_0\ldots k_{m-1}}-(\Delta_R^{m-1}y^{k_m})_{i_0\ldots i_{m-2}i_m}^{k_0\ldots k_{m-1}})^{-1}.
\end{aligned}
\end{equation}

b) {\it We have}
\begin{equation}
\begin{aligned}{}
&(\Delta_L^0 g)_i^k=(y_i^k)^{-1}g^k,\\
&(\Delta_L^1 g)_{ii'}^{kk'}=((\Delta_L^0 y_{i'})_i^k-(\Delta_L^0 y_{i'})_i^{k'})^{-1}((\Delta_L^0 g)_i^k-(\Delta_L^0 g)_i^{k'}),\\
&(\Delta_L^m g)_{i_0\ldots i_m}^{k_0\ldots k_m}=((\Delta_L^{m-1}y_{i_m})_{i_0\ldots i_{m-1}}^{k_0\ldots k_{m-1}}-(\Delta_L^{m-1}y_{i_m})_{i_0\ldots i_{m-1}}^{k_0\ldots k_{m-2}k_m})^{-1}\\
&\times((\Delta_L^{m-1}g)_{i_0\ldots i_{m-1}}^{k_0\ldots k_{m-1}}-(\Delta_L^{m-1}g)_{i_0\ldots i_{m-1}}^{k_0\ldots k_{m-2}k_m}).
\end{aligned}
\end{equation}

{\it Proof.} Let us prove part (b). The proof of (a) is similar.

Let us fix the sequence $i_0,i_1,\ldots$.
Denote by $w^m$ and $w^{m-1}$ the biorthogonal basis vectors corresponding to the sequence $k_0$, $\ldots$, $k_{m-1}$, $k_m$, and by $\widetilde w^{m-1}$ the basis vector corresponding to the sequence 
$k_0$, $\ldots$, $k_{m-2}$, $ k_m$.
We have
$$
\begin{aligned}{}
&(w^{m-1},p_{i_j})=(\widetilde w^{m-1},p_{i_j})=0,\ \  0\le j\le m-2,\\
&(w^{m-1},p_{i_{m-1}})=(\widetilde w^{m-1},p_{i_{m-1}})=1.
\end{aligned}
$$
This implies that $w^{m-1}-\widetilde w^{m-1}$ is a linear combination of $v^{k_0}$, $\ldots$, $v^{k_m}$ satisfying the properties 
$$
(w^{m-1}-\widetilde w^{m-1}, p_{i_j})=0,\ \ 0\le j\le m-1.
$$
Therefore, we have $w^m=c(w^{m-1}-\widetilde w^{m-1})$. Pairing this equality with $p_{i_m}$, we find $c=((w^{m-1},p_{i_m})-(\widetilde w^{m-1},p_{i_m}))^{-1}$ and hence
$$
(w^m,\g)=((w^{m-1},p_{i_m})-(\widetilde w^{m-1},p_{i_m}))^{-1}((w^{m-1},\g)-(\widetilde w^{m-1},\g)).
$$ 
In notations (8), this is formula (10). \qed
\medskip

{\it Explicit formulas for difference derivatives.}

{\bf Theorem 3.} a) {\it $(\Delta_R^m f)_{i_0\ldots i_m}^{k_0\ldots k_m}$ is the unique right linear combination of $f_{i_0}$, $\ldots$, $f_{i_m}$ with the properties}
\begin{equation}
(\Delta_R^m y^{k_l})_{i_0\ldots i_m}^{k_0\ldots k_m}=0,\ \ 0\le l\le m-1;\ \ (\Delta_R^m y^{k_m})_{i_0\ldots i_m}^{k_0\ldots k_m}=1. 
\end{equation}
b) {\it $(\Delta_L^m g)_{i_0\ldots i_m}^{k_0\ldots k_m}$ is the unique left linear combination of $g^{k_0}$, $\ldots$, $g^{k_m}$ with the properties}
\begin{equation}
(\Delta_L^m y_{i_j})_{i_0\ldots i_m}^{k_0\ldots k_m}=0,\ \ 0\le j\le m-1;\ \ (\Delta_L^m y_{i_m})_{i_0\ldots i_m}^{k_0\ldots k_m}=1. 
\end{equation}

{\it Proof} follows from equalities (7, 8).\qed
\medskip

{\bf Theorem 4.} {\it We have}

a) $(\Delta_R^m f)_{i_0\ldots i_m}^{k_0\ldots k_m}=\sum\limits_{j=0}^m f_{i_j}(|D_{i_0\ldots i_m}^{k_0\ldots k_m}|_j^m)^{-1}$;

b) $(\Delta_L^m g)_{i_0\ldots i_m}^{k_0\ldots k_m}=\sum\limits_{l=0}^m(|D^{k_0\ldots k_m}_{i_0\ldots i_m}|_m^l)^{-1}g^{k_l}$.
\medskip

{\bf Corollary.} a) {\it $(\Delta_R^m f)_{i_0\ldots i_m}^{k_0\ldots k_m}$ is symmetric with respect to $i_0,\ldots,i_m$ and with respect to 
$k_0,\ldots,k_{m-1}$.}

b) {\it  $(\Delta_L^m g)_{i_0\ldots i_m}^{k_0\ldots k_m}$ is symmetric with respect to $k_0,\ldots,k_m$ and with respect to $i_0,\ldots,i_{m-1}$.}

\subsection{The case $i_m=k_m=m$} Let us fix 
$$
(i_0,i_1,\ldots)=(k_0,k_1,\ldots)=(0,1,\ldots). 
$$
Denote $D_{0\ldots n}^{0\ldots n}=D_n$, $A_{0\ldots n}^{0\ldots n}=A_n$,
$C_{0\ldots n}^{0\ldots n}=C_n$;
$$
(\Delta_R^m f)_{0\ldots m}^{0\ldots m}=\Delta_R^m f;\ \  (\Delta_L^m g)_{0\ldots m}^{0\ldots m}=\Delta_L^m g.
$$

By Theorem 3, 
\begin{equation}
\begin{aligned}{}
\Delta_R^m f&=(\f,q_m)=\sum\limits_{i=0}^m f_i a_m^i,\\
\Delta_L^m g&=(w^m,\g)=\sum\limits_{k=0}^m c_k^m g^k
\end{aligned}
\end{equation}
are the unique linear combinations of respectively $f_0,\ldots,f_m$ and $g^0$, $\ldots$, $g^m$ with the properties
\begin{equation}
\begin{aligned}{}
&\Delta_R^m y^k=0,\ \ 0\le k\le m-1,\ \ \Delta_R^m y^m=1;\\
&\Delta_L^m y_i=0,\ \  0\le i\le m-1,\ \ \Delta_L^m y_m=1.
\end{aligned}
\end{equation}
By Theorem 4, we have
\begin{equation}
\begin{aligned}{}
\Delta_R^m f&=\sum\limits_{i=0}^m f_i(|D_m|_i^m)^{-1},\\
\Delta_L^m g&=\sum\limits_{k=0}^m(|D_m|_m^k)^{-1}g^k.
\end{aligned}
\end{equation}
\medskip

{\bf Theorem 5.}
a) {\it  $A_n$ is the unique upper triangular matrix such that the matrix $D_n A_n$ is strictly lower triangular 
\emph(i.~e. it is a lower triangular matrix with the units on the main diagonal\emph{);}}

 b) {\it  $C_n$ is the unique lower triangular matrix such that the matrix $C_n D_n$ is strictly upper triangular 
\emph(i.~e. it is an upper triangular matrix with the units on the main diagonal\emph).}

{\it Proof} follows from equalities (13, 14).\qed

\subsection{Generalized Lagrange interpolation formula} Put
\begin{equation}
(\f,\g)_n=\sum\limits_{i,k=0}^n (\f,p_i)(z_k^i)_n(v^k,\g)=\sum\limits_{i,k=0}^n f_i(z_k^i)_n g^k,
\end{equation}
where $(z_k^i)_n$ is the inverse matrix to the matrix $D_n$,
\begin{equation}
(z_k^i)_n=(|D_n|_i^k)^{-1}.
\end{equation}
We have
\begin{equation}{}
\begin{aligned}{}
&(\f,p_i)_n=f_i,\ \ 0\le i\le n;\\
&(v^k,\g)_n=g^k,\ \ 0\le k\le n.
\end{aligned}
\end{equation}
Hence, 
\begin{equation}
\begin{aligned}{}
&(\f,\g)=(\f,\g)_n+R_n(\f,\g),\\
&R_n(\f,p_i)=R_n(v^k,\g)=0,\ \ 0\le i,k\le n.
\end{aligned}
\end{equation}

\subsection{Formula for the inverse matrix}
{\bf Theorem 6.}
\begin{equation}
(z_k^i)_n=\sum\limits_{m=\max(i,k)}^n (|D_m|_i^m)^{-1}|D_m|_m^m(|D_m|_m^k)^{-1}.
\end{equation}

{\it Proof} follows from the equalities, for $i,k<m$,
$$
\begin{aligned}{}
&(z_k^i)_m-(|D_m|_i^m)^{-1}|D_m|_m^m(|D_m|_m^k)^{-1}\\
&=(z_k^i)_m-(z_m^i)_m(z_m^m)_m^{-1}(z_k^m)_m\\
&=\left|\begin{array}{cc}(z_k^i)_m&(z_m^i)_m\\(z_k^m)_m&(z_m^m)_m\end{array}\right|_1^1=(|D_{m-1}|_i^k)^{-1}\text{ (see [3])}\\
&=(z_k^i)_{m-1}.\ \ \qed
\end{aligned}
$$

\subsection{Generalized Newton formula} {\bf Theorem 7.} 
\begin{equation}
\begin{aligned}{}
&(\f,\g)=\sum\limits_{m=0}^n\Delta_R^m f|D_m|_m^m\Delta_L^m g+R_n(\f,\g),\\
&R_n(\f,p_i)=R_n(v^k,\g)=0,\ \ 0\le i,k\le n.
\end{aligned}
\end{equation}

{\it Proof.} We shall prove the equality
\begin{equation}
(\f,\g)_n=\sum\limits_{m=0}^n\Delta_R^m f|D_m|_m^m\Delta_L^m g
\end{equation}
and use (19).

{\it First proof of} (22).
$$
\begin{aligned}{}
&(\f,\g)_n=\sum\limits_{i,k=0}^n f_i(z_k^i)_n g^k\\
&\stackrel{\rm Th. 6}{=}\sum\limits_{m=0}^n\sum\limits_{i,k=0}^m f_i(|D_m|_i^m)^{-1}|D_m|_m^m(|D_m|_m^k)^{-1}g^k\\
&\stackrel{\rm (15)}{=}\sum\limits_{m=0}^n\Delta_R^m f|D_m|_m^m\Delta_L^m g. \qed
\end{aligned}
$$

{\it Second proof.} By Theorem 5 we have $D_n^{-1}=A_n B_n C_n$, where $A_nB_n$ is a strictly upper triangular matrix and $B_nC_n$ is a strictly lower
triangular matrix. This implies that $B_n=\diag(b_m^m)$ is a diagonal matrix. By formula (15), we have $a_m^m=c_m^m=(b_m^m)^{-1}=(|D_m|_m^m)^{-1}$. This implies 
formula (22).\qed

{\it Third proof.} Since the bases $w^m$, $q_m$ are biorthogonal, we have
$$
\begin{aligned}{}
(\f,\g)_n&=\sum\limits_{m=0}^n(\f,q_m)(w^m,q_m)^{-1}(w^m,\g)\\
&=\sum\limits_{m=0}^n\Delta_R^m f(w^m,q_m)^{-1}\Delta_L^m g.
\end{aligned}
$$
By (13, 15), we have $(w^m,q_m)=(|D_m|_m^m)^{-1}$.\qed

\section{Examples and applications}

\subsection{Example 1: the Newton formula} Let $V$ be a space of analytical functions $\f=f(x)$, let $v^k(x)=x^k$, let $x_0, x_1,\ldots$ be pairwise different numbers,  
let $(\f,p_i)=f(x_i)$, and let $(\f,\g)=f(x)$. Then we have $y_i^k=x_i^k$ (the $k$-th power of $x_i$), $f_i=f(x_i)$, $g^k=x^k$.

By (13, 14), $\Delta_L^m g=P_m(x)$ is the unique degree $m$ polynomial in $x$ with the properties
\begin{equation}
\begin{aligned}{}
&P_m(x_0)=P_m(x_1)=\ldots=P_m(x_{m-1})=0,\\ 
&P_m(x_m)=1.
\end{aligned}
\end{equation}
Therefore, we have
\begin{equation}
\Delta_L^m g=\frac{(x-x_0)\ldots(x-x_{m-1})}{(x_m-x_0)\ldots(x_m-x_{m-1})}.
\end{equation}
Further, by (9), $\Delta_R^m f=\Delta^m f(x_0,\ldots,x_m)$ is the $m$-th difference derivative (1).
More generally, $(\Delta_R^m  y^k)_{0\ldots m}^{k_0\ldots k_m}$ is a ratio of Schur symmetric functions in $x_0,\ldots,x_m$. 

Also, we have 
\begin{equation}
|D_m|_m^m=(x_m-x_0)\ldots(x_m-x_{m-1});
\end{equation}
and $(\f,\g)_n=S_n(x)$ is the unique degree $n$ polynomial in $x$ with the properties 
\begin{equation}
S_n(x_i)=f_i,\ \ 0\le i\le n,
\end{equation}
i.~e. it is the Lagrange interpolation polynomial. Hence, in this case the formula of Theorem 7 is the Newton formula (2).

\subsection{Example 2: the Taylor formula} Let $V$ be a space of smooth functions $\f=f(x)$, let $v^k(x)=x^k$, let $x_0$ be a number, let
$(\f,p_i)=f^{(i)}(x_0)$, and let $(\f,\g)=f(x)$. Then we have $y_i^k=\frac{d^i x^k}{dx^i}(x_0)$ (the $i$-th derivative of $x^k$ at $x=x_0$), $f_i=f^{(i)}(x_0)$, $g^k=x^k$.

By (13, 14), $\Delta_L^m g=P_m(x)$ is the unique degree $m$ polynomial in $x$ with the properties
\begin{equation}
\begin{aligned}{}
&P_m(x_0)=P_m'(x_0)=\ldots=P_m^{(m-1)}(x_0)=0,\\ 
&P_m^{(m)}(x_0)=1.
\end{aligned}
\end{equation}
Therefore, 
\begin{equation}
\Delta_L^m g=(x-x_0)^m/m!.
\end{equation}
Further, by (14) we have $\Delta_R^m f=f^{(m)}(x_0)/m!$;
\begin{equation} 
|D_m|_m^m=m!;
\end{equation} 
and $(\f,\g)_n=S_n(x)$ is the unique degree $n$ polynomial in $x$ with the properties
\begin{equation}
S_n^{(i)}(x_0)=f^{(i)}(x_0),\ \ 0\le i\le n,
\end{equation}
i.~e. it is the Taylor polynomial. Hence, in this case the formula of Theorem 7 is the Taylor formula (4).

As $x_1,\ldots,x_m\to x_0$, the polynomial $|D_m|_m^m\Delta_L^m g$ given by (24, 25) tends to the polynomial given by (28, 29).
This yields an explanation of the limit formula (3).

\subsection{Example 3: the orthogonalization process}
Let $V=V'$ be a Hilbert space with the scalar product $(,)$, and let $v^k=p_k$ be a basis in $V$. 
There exists the unique nondegenerate lower triangular matrix $C=(c_k^m)$ such 
that the basis 
\begin{equation}
w^m=\sum\limits_{k=0}^m c_k^m v^k
\end{equation}
in $V$ is orthogonal and $(w^m,v^m)=1$. The matrix $C$ and the basis $w^m$ are obtained by applying the orthogonalization process to the basis $v^k$. Then 
the formula of Theorem 7 amounts to the formula
\begin{equation}
(\f,\g)=\sum\limits_m (\f,w^m)(w_m,w_m)^{-1}(w_m,\g).
\end{equation}
This case and its various applications to analysis are studied in detail in the book [1].

\end{document}